\documentclass[12pt, a4paper]{amsart}

\usepackage[all]{xy}
\usepackage{amsthm}
\usepackage{amsmath}
\usepackage{amssymb}
\usepackage{amsfonts}
\usepackage{graphicx}
\usepackage[mathscr]{eucal}
\usepackage{enumerate}
\usepackage[latin1]{inputenc}

\newtheorem{thm}{Theorem}[section]
\newtheorem{prop}[thm]{Proposition}

\newtheorem{cor}[thm]{Corollary}
\newtheorem{que}[thm]{Question}
\newtheorem{prob}[thm]{Problem}
\newtheorem{exmp}[thm]{Example}
\newtheorem{lem}[thm]{Lemma}
\newtheorem{rem}[thm]{Remark}

\newtheorem{conj}{Conjeture}

\newenvironment{pf}{{\noindent\sc Proof\;}}{\qed\\}

\def\length{\mathrm{\underline{length}}\hspace{0.1cm}}
\def\len{\mathrm{length}\hspace{0.1cm}}

\def\Sup{\mathrm{PSupp}}

 \DeclareMathOperator{\Hom}{Hom}
\def\Coker{\mathrm{Coker\hspace{0.1cm}}}
\def\Cohom{\mathrm{Cohom}}

\def\Ext{\mathrm{Ext}}
\def\Hom{\mathrm{Hom}}

\def\dim{\mathrm{dim}}

\def\Soc{\mathrm{soc\hspace{0.1cm}}}

\def\Ima{\mathrm{Im\hspace{0.1cm}}}

\def\T{\mathcal{T}}
\def\K{\mathcal{K}}
\def\M{\mathcal{M}}

\title{Localization in Tame and Wild Coalgebras}
\author{Pascual Jara}
\author{Luis M. Merino}
\author{Gabriel Navarro}
\address{Department of Algebra,
University of Granada\\
Avda. Fuentenueva s/n, E-18071, Granada, Spain} \email{pjara@ugr.es,
lmerino@ugr.es, gnavarro@ugr.es}
\thanks{Research partially supported by DGES BMF2001-2823 and FQM-266
(Junta de Andaluc{\'\i}a Research Group). The third author is
partially supported by Spanish MEC grant BES-2002-2403.}

\date{}

\begin{document}
\maketitle

\begin{abstract}
We apply the theory of localization for tame and wild coalgebras in
order
 to prove the following theorem: "Let $Q$ be an acyclic quiver. Then any tame admissible subcoalgebra
of $KQ$ is the path coalgebra of a quiver with relations".
\end{abstract}

\section{Introduction}

The quiver-theoretical ideas developed by Gabriel  and his school
during the seventies have been the origin of many advances in
Representation Theory of Algebras for years. Moreover, many of the
present developments of the theory use, up to some extent, these
techniques and results. Among these tools, it is mostly accepted
that the main one is the famous Gabriel theorem which relates any
finite dimensional algebra (over an algebraically closed field) to a
nice quotient of a path algebra (see for instance \cite{simsonbook1}
or \cite{auslander}). The key-point of such method lies in the
achievement of a description of all finitely generated modules by
means of finite dimensional representations of the quiver.

More recently, some authors have tried to get rid of the imposed
finiteness conditions by the use of coalgebras and their category of
comodules, see \cite{chin}, \cite{justus}, \cite{simson1},
\cite{simson2} and \cite{woodcock}. Coalgebras are realized, because
of the freedom on choosing their dimension, as an intermediate step
between the representation theory of finite dimensional and infinite
dimensional algebras. Also, because of their locally finite nature,
coalgebras are a good candidate for extending many techniques and
results stated for finite dimensional algebras. In \cite{jmn}, the
authors intent to prove that every basic coalgebra, over an
algebraically closed field $K$, is the path coalgebra $C(Q,\Omega)$
of a quiver $Q$ with a set of relations $\Omega$, in the sense of
Simson \cite{simson1}. That is, an analogue for coalgebras of
Gabriel's theorem. The result is proven in \cite{simson2} for the
family of coalgebras $C$ such that the Gabriel quiver $Q_C$ of $C$
is intervally finite. Unfortunately, that proof does not hold for
arbitrary coalgebras, as a class of counterexamples given in
\cite{jmn} shows. Moreover, there it is proven a criterion allowing
us to decide whether or not a coalgebra is of that kind (see Theorem
\ref{criterion} below).

Nevertheless, it is worth noting that for all counterexamples found
in \cite{jmn}, the category of finitely generated comodules has very
bad properties. This bad behaviour is similar to the notion of
wildness given for the finitely generated module categories of some
finite dimensional algebras, meaning that this category is so big
that it contains (via an exact representation embedding) the
category of all finite dimensional representations of the
noncommutative polynomial algebra $K\langle x,y\rangle$. As it is
well-known, the category of finite dimensional modules over
$K\langle x,y\rangle$ contains (again via an exact representation
embedding) the category of all finitely generated representations
for any other finite dimensional algebra, and thus it is not
realistic aiming to give an explicit description of this category
(or, by extension, of any wild algebra). The counterpart to the
notion of wild algebra is the one of tameness, a tame algebra being
one whose indecomposable modules of finite dimension are
parametrised by a finite number of one-parameter families for each
dimension vector. A classical result in representation theory of
algebras (the Tame-Wild Dichotomy, see \cite{drozd}) states that any
finite dimensional algebra over an algebraically closed field is
either of tame module type or of wild module type. We refer the
reader to \cite{simsonbluebook} for basic definitions and properties
about module type of algebras.

Analogous concepts were defined by Simson in \cite{simson1} for
coalgebras, see Section~\ref{deftame}. In \cite{simson2}, it was
proven a weak version of the Tame-Wild Dichotomy that goes as
follows: over an algebraically closed field $K$, if $C$ is
$K$-coalgebra of tame comodule type, then $C$ is not of wild
comodule type. The full version remains open:

\begin{conj} Any $K$-coalgebra, over an algebraically closed field
$K$, is either of tame comodule type, or of wild comodule type, and
these types are mutually exclusive.
\end{conj}

As, according to Theorem \ref{criterion} below, the coalgebras which
are not a path coalgebra of a quiver with relation are close to be
wild, we may reformulate the problem stated in \cite[Section
8]{simson1} as follows:

\begin{que}\label{que1} Assume that $K$ is an algebraically closed field. Is any
basic coalgebra $C$ of tame comodule type isomorphic to the path
coalgebra $C(Q,\Omega)$ of a quiver $Q$ with a set of relations
$\Omega$?
\end{que}

One of the main results of this paper (Corollary \ref{maincor}, as a
consequence of Theorem \ref{thm}) shows that the answer to Question
\ref{que1} is affirmative, under the further assumption that the
quiver $Q$ is acyclic, that is, $Q$ has no oriented cycles. This
extends \cite[Theorem 3.14 (c)]{simson2} from the case in which $Q$
is intervally finite and $C'\subseteq KQ$ is an arbitrary admissible
subcoalgebra to the case in which $Q$ is acyclic and $C'\subseteq
KQ$ is tame.

In the proof of Theorem \ref{thm} and Corollary \ref{maincor} we use
the theory of localization introduced by Gabriel \cite{gabriel}, see
also \cite{blas2} and \cite{popescu}. For that reason, in Section
\ref{sectpointed}, we develop that theory for arbitrary pointed
coalgebras (subcoalgebras of path algebras) extending the results
stated for path coalgebras (hereditary pointed coalgebras) in
\cite[Section 3]{jmnr}. Then we show in Section \ref{sectiontame}
and Section \ref{sectionwild} that, under suitable assumptions, the
tameness and wildness of a coalgebra $C$ is preserved by
localization. Lastly, in Section \ref{sectiongabriel}, we prove the
above-mentioned result: for any algebraically closed field $K$ and
any acyclic quiver $Q$, every tame admissible subcoalgebra $C'$ of
the path coalgebra $KQ$ of the quiver $Q$ is isomorphic to the path
coalgebra $C(Q,\Omega)$ of quiver with relations $(Q,\Omega)$.

\section{Preliminaries}

Throughout this paper $K$ will be an algebraically closed field. By
a quiver, $Q$, we mean a quadruple $(Q_0,Q_1,h,s)$ where $Q_0$ is
the set of vertices (points), $Q_1$ is the set of arrows and for
each arrow $\alpha\in Q_1$, the vertices $h(\alpha)$ and $s(\alpha)$
are the source (or start point) and the sink (or end point) of
$\alpha$, respectively (see \cite{simsonbook1} and
\cite{auslander}). If $i$ and $j$ are vertices in $Q$, an (oriented)
path in $Q$ of length $m$ from $i$ to $j$ is a formal composition of
arrows
$$p=\alpha_m\cdots \alpha_2 \alpha_1$$ where $h(\alpha_1)=i$,
$s(\alpha_m)=j$ and $s(\alpha_{k-1})=h(\alpha_{k})$, for $k=2,
\ldots , m$. To any vertex $i\in Q_0$ we attach a trivial path of
length 0, say $e_i$, starting and ending at $i$ such that $\alpha
e_i=\alpha$ (resp. $e_i\beta=\beta$) for any arrow $\alpha$ (resp.
$\beta$) with $h(\alpha)=i$ (resp. $s(\beta)=i$). We identify the
set of vertices and the set of trivial paths. An (oriented) cycle is
a path in $Q$ which starts and ends at the same vertex. $Q$ is said
to be acyclic if there is no oriented cycle in $Q$

Let $KQ$ be the $K$-vector space generated by the set of all paths
in $Q$. Then $KQ$ can be endowed with the structure of a (non
necessarily unitary) $K$-algebra with multiplication induced by
concatenation of paths, that is,
$$(\alpha_m\cdots \alpha_2\alpha_1)(\beta_n\cdots
\beta_2\beta_1)=\left\{%
\begin{array}{ll}
    \alpha_m\cdots \alpha_2\alpha_1\beta_n\cdots\beta_2\beta_1 & \hspace{0.2cm}\text{if $s(\beta_n)=h(\alpha_1),$} \\
    0 & \hspace{0.2cm}\text{otherwise;} \\
\end{array}%
\right.$$ $KQ$ is the path algebra of the quiver $Q$. The algebra
$KQ$ can be graded by
$$KQ=KQ_0\oplus KQ_1\oplus \cdots \oplus KQ_m\oplus \cdots,$$
where $Q_m$ is the set of all paths of length $m$. An ideal
$\Omega\subseteq KQ$ is called an ideal of relations or a relation
ideal if $\Omega\subseteq KQ_2\oplus KQ_3\oplus\cdots =KQ_{\geq 2}$.
By a quiver with relations we mean a pair $(Q,\Omega)$, where $Q$ is
a quiver and $\Omega$ a relation ideal of $KQ$. For more details and
basic facts from representation theory of algebras the reader is
referred to \cite{simsonbook1} and \cite{auslander}.

Following \cite{woodcock}, the path algebra $KQ$ can be viewed as a
graded $K$-coalgebra with comultiplication induced by the
decomposition of paths, that is, if $p=\alpha_m\cdots \alpha_1$ is a
path from the vertex $i$ to the vertex $j$, then
$$\Delta(p)=e_j\otimes p+p\otimes e_i+\sum_{i=1}^{m-1}\alpha_m\cdots
 \alpha_{i+1}\otimes
\alpha_i\cdots \alpha_1=  \sum_{\eta \tau=p} \eta\otimes\tau$$ and
for a trivial path, $e_i$, we have $\Delta(e_i)=e_i\otimes e_i$. The
counit of $KQ$ is defined by the formula
$$\epsilon(\alpha)=
\left\{%
\begin{array}{ll}
    1 & \text{if $\alpha\in Q_0$}, \\
    0 & \text{if $\alpha$ is a path of length $\geq 1$}. \\
\end{array}%
\right.
$$
The coalgebra $(KQ,\Delta,\epsilon)$ (shortly $KQ$) is called the
path coalgebra of the quiver $Q$. A subcoalgebra $C$ of a path
coalgebra $KQ$ of a quiver $Q$ is said to be admissible if it
contains the subcoalgebra of $KQ$ generated by all vertices and all
arrows.

Let $(Q,\Omega)$ be a quiver with relations. Following
\cite{simson1}, the path coalgebra of $(Q,\Omega)$ is defined by the
subspace of $KQ$,
$$C(Q,\Omega)=\{\text{$a\in KQ$ $\mid$ $\langle a,\Omega
\rangle =0$} \}$$ where $\langle -,- \rangle:KQ\times KQ
\longrightarrow K$ is the bilinear map defined by $\langle
v,w\rangle =\delta_{v,w}$ (the Kronecker delta) for any two paths
$v$ and $w$ in $Q$.

Let us recall from \cite{simson2} the following weak version of
Gabriel theorem for coalgebras. A quiver is called intervally finite
if there are finitely many paths between any two vertices.

\begin{prop}
Let $Q$ be an intervally finite quiver. Then the map $\Omega
\longmapsto C(Q,\Omega )$ defines a bijection between the set of all
relation ideals $\Omega$ of $KQ$ and the set of all admissible
subcoalgebras of $KQ$.
\end{prop}

As mentioned in the previous section, the following  criterion is
proven in \cite{jmn}:

\begin{thm}\label{criterion}
Let $C$ be an admissible subcoalgebra of a path coalgebra $KQ$. Then
$C$ is not the path coalgebra of a quiver with relations if and only
if there exist an infinite number of different paths
$\{\gamma_{i}\}_{i\in \mathbb{N}}$ in $Q$ such that:
\begin{enumerate}[$(a)$]
\item All of them have common source and common sink.
\item None of them is in $C$.
\item There exist elements $a_j^n\in K$ for all $j,n\in
\mathbb{N}$ such that the set $\{\gamma_n+\sum_{j>n} a^n_j
\gamma_j\}_{n\in \mathbb{N}}$ is contained in $C$.
\end{enumerate}
\end{thm}

\begin{rem} Following the construction of the elements of condition
$(c)$ exposed in \cite{jmn}, the sum $\sum_{j>n} a^n_j \gamma_j$ is
finite for any $n\in \mathbb{N}$.
\end{rem}

Given a $K$-coalgebra $C$, we denote by $\M^C_f$, $\M^C_{qf}$ and
$\M^C$ the category of finite dimensional, quasi-finite and all
right $C$-comodules, respectively. Let now  $\T$ be a dense
subcategory of the category of right $C$-comodules. Following
\cite{gabriel} and \cite{popescu}, $\T$ is said to be localizing if
the quotient functor $T:\M^C\rightarrow \M^C/\T$ has a right adjoint
functor, $S$, called the section functor.  If the section functor is
exact, $\T$ is called perfect localizing. Dually, see \cite{blas2},
$\T$ is said to be colocalizing if $T$ has a left adjoint functor,
$H$, called the colocalizing functor. $\T$ is said to be perfect
colocalizing if the colocalizing functor is exact.

Let us list some properties of the localizing functors (see
\cite{gabriel} and \cite{blas2}).

\begin{lem} Let $\T$ a dense subcategory
of the category of right comodules $\M^C$ over a coalgebra $C$. The
following statements hold:
\begin{enumerate}[$(a)$]
\item $T$ is exact.
\item If $\T$ is localizing, then the section functor $S$ is left exact and
the equivalence $TS\simeq 1_{\M^C/\T}$ holds.
\item If $\T$ is colocalizing, then the colocalizing functor $H$ is right exact and
the equivalence $TH\simeq 1_{\M^C/\T}$ holds.
\end{enumerate}
\end{lem}

From the general theory of localization in Grothendieck categories,
it is well-known that there exists a one-to-one correspondence
between localizing subcategories of $\M^C$ and sets of
indecomposable injective right $C$-comodules, and, as a consequence,
sets of simple right $C$-comodules. In \cite{cuadra}, \cite{jmnr}
and \cite{woodcock}, localizing subcategories are described by means
of idempotents in the dual algebra $C^*$. In particular, it is
proved that the quotient category is the category of right comodules
over the coalgebra $eCe$, where $e$ is the idempotent associated to
the localizing subcategory. The coalgebra structure of $eCe$ is
given by
$$\Delta_{eCe} (exe)=\displaystyle \sum_{(x)} ex_{(1)}e\otimes
ex_{(2)}e \hspace{0.4cm}\text{ and }\hspace{0.4cm}
\epsilon_{eCe}(exe)=\epsilon_C(x)$$ for any $x\in C$, where
$\Delta_C(x)= \sum_{(x)} x_{(1)} \otimes x_{(2)}$ using the
sigma-notation of \cite{sweedler}. Throughout we denote by $\T_e$
the localizing subcategory associated to the idempotent $e$. For
completeness, we remind from \cite{cuadra} (see also \cite{jmnr})
the following description of the localizing functors. We recall
that, given an idempotent $e\in C^*$, for each right $C$-comodule
$M$, the vector space $eM$ is endowed with a structure of right
$eCe$-comodule given by
$$\rho_{eM}(ex)=\sum_{(x)} ex_{(1)}\otimes
ex_{(0)} e$$ where $\rho_{M}(x)=\sum_{(x)} x_{(1)}\otimes x_{(0)}$
using the sigma-notation of \cite{sweedler}.

\begin{lem} Let $C$ be a coalgebra and $e$ be an idempotent in
$C^*$. Then the following statements hold:
\begin{enumerate}[$(a)$]
\item The quotient functor $T:\M^C\rightarrow \M^{eCe}$ is naturally
equivalent to the functor $e(-)$. $T$ is also equivalent to the
cotensor functor $-\square_{C} eC$.
\item The section functor $S:\M^{eCe} \rightarrow \M^C$ is naturally
equivalent to the cotensor functor $-\square_{eCe} Ce$.
\item If $\T_e$ is a colocalizing subcategory of $\M^C$, the
colocalizing functor $H:\M^{eCe}\rightarrow \M^C$ is naturally
equivalent to the functor $\Cohom_{eCe}(eC,-)$.
\end{enumerate}
\end{lem}

Dealing with quivers, the localization in categories of comodules
over path coalgebras is described in detail in \cite[Section
3]{jmnr}. If $Q=(Q_0,Q_1)$ is a quiver and $C$ is an admissible
subcoalgebra of the path coalgebra $KQ$, localizing subcategories of
$\M^{C}$ are given by subsets $X$ of the set of vertices $Q_0$. We
denote by $e_X$ the idempotent of $C^*$ related to the set $X$, and
the subcategory $\T_X:=\T_{e_X}$ of $\M^C$ is called the localizing
subcategory associated to $X$.

For the convenience of the reader we summarize the functors obtained
in the situation of (co)localization by means of the diagrams

$$ \xymatrix@C=50pt{ \mathcal{M}^{C}
\ar@<0.75ex>[rr]^-{T=e(-)=-\square_{C}eC} & &
\ar@<0.75ex>[ll]^-{S=-\square_{eCe}Ce} \mathcal{M}^{eCe}}
\hspace{0.2cm}\text{and}\hspace{0.2cm}
 \xymatrix@C=50pt{ \mathcal{M}^{C}
\ar@<0.75ex>[rr]^-{T=e(-)=-\square_{C}eC} & &
\ar@<0.75ex>[ll]^-{H=\Cohom_{eCe}(eC,-)} \mathcal{M}^{eCe}}.$$

\section{Localization in pointed coalgebras}\label{sectpointed}

This section is devoted to develop the theory of localization for
pointed coalgebras. In particular, we generalize the results given
in \cite[Section 3]{jmnr} for path coalgebras to an arbitrary
subcoalgebra of a path coalgebra. We recall that, over an
algebraically closed field, every coalgebra is Morita-Takeuchi
equivalent to a pointed coalgebra.

For the convenience of the reader we introduce the following
notation:

Let $Q$ be a quiver and $p=\alpha_n\alpha_{n-1}\cdots\alpha_1$ be a
non-trivial path in $Q$. We denote by $I_p$ the set of vertices
$\{h(\alpha_1), s(\alpha_1), s(\alpha_2),\ldots ,s(\alpha_n)\}$.
Given a subset of vertices $X\subseteq Q_0$, we say that $p$ is a
\emph{cell} in $Q$ relative to $X$ (shortly a cell) if $I_p\cap
X=\{h(p), s(p)\}$ and $s(\alpha_i)\notin X$ for all $i=1,\ldots
,n-1$. Given $x, y\in X$, we denote by $\mathcal{C}ell^Q_X(x,y)$ the
set of all cells from $x$ to $y$. We denote the set of all cells in
$Q$ relative to $X$ by $\mathcal{C}ell^Q_X$.

\begin{lem} Let $Q$ be a quiver and $X\subseteq Q_0$ be a subset of
vertices. Given a path $p$ in $Q$ such that $h(p)$ and $s(p)$ are in
$X$, then $p$ has a unique decomposition $p=q_r\cdots q_1$, where
each $q_i$ is a cell in $Q$ relative to $X$.
\end{lem}
\begin{pf}
It is straightforward.
\end{pf}

Let $p$ be a non-trivial path in $Q$ which starts and ends at
vertices in $X\subseteq Q_0$. We call the \emph{cellular
decomposition}\index{cellular decomposition} of $p$ relative to $X$
to the decomposition given in the above lemma.

Following \cite{woodcock}, every pointed coalgebra $C$ is isomorphic
to a subcoalgebra of the path coalgebra $KQ$ of the Gabriel quiver
$Q_C$ of $C$. Furthermore, $C$ is an admissible subcoalgebra and
then, by \cite[Corollary 3.5]{jmn}, it is a relation subcoalgebra in
the sense of Simson \cite{simson2}, that is, $C$ has a
decomposition, as vector space, $C=\bigoplus_{a,b\in Q_0} C_{ab}$,
where $C_{ab}=C\cap KQ(a,b)$ and $Q(a,b)$ is the set of all paths in
$Q$ from $a$ to $b$. Hence, $C$ has a basis in which every basic
element is a linear combination of paths with common source and
common sink.

Fix an idempotent $e\in C^*$. Then there exists a subset $X$ of
$Q_0$ such that the localized coalgebra $eCe$ has a decomposition
$eCe=\bigoplus_{a,b\in X} C_{ab}$, that is, the elements of $eCe$
are linear combinations of paths with source and sink at vertices in
$X$. It follows that $eCe$ is a pointed coalgebra so, by
\cite{woodcock} again, there exists a quiver $Q^e$ such that $eCe$
is an admissible subcoalgebra of $KQ^e$. The quiver $Q^e$ is
described as follows:

\emph{Vertices}. We know that $Q_0$ equals the set of group-like
elements
 $\mathcal{G}(C)$ of $C$, therefore
 $(Q^e)_0=\mathcal{G}(eCe)=e\mathcal{G}(C)e=e Q_0 e=X$.

\emph{Arrows}. Let $x$ and $y$ be vertices in $X$. An element $p\in
eCe$ is
  a non-trivial $(x,y)$-primitive element if and only if $p\notin
  KX$ and $\Delta_{eCe}(p)=y\otimes p+p\otimes x$. Without loss of generality we may assume that
  $p=\sum_{i=1}^{n} \lambda_i p_i$ is an element in $eCe$ such that each path $p_i$ is not
  trivial, and
   $p_i=\alpha^i_{n_i}\cdots\alpha_2^i\alpha^i_1$ and
  $p_i=q^i_{r_i}\cdots q_1^i$ are the decomposition of $p_i$ in
  arrows of $Q$ and the cellular decomposition of $p_i$ relative to
  $X$, respectively, for all $i=1,\ldots, n$.  Then
$$\Delta_C(p)=\sum_{i=1}^n \lambda_i p_i\otimes h(p_i)+
\sum_{i=1}^n \lambda_i s(p_i)\otimes p_i+\sum_{i=1}^n \lambda_i
\sum_{j=2}^{n_i} \alpha^i_{n_i}\cdots
\alpha^i_j\otimes\alpha^i_{j-1}\cdots\alpha^i_1
$$
 and therefore,
\begin{multline*}\Delta_{eCe}(p)=\sum_{i=1}^n \lambda_i (e\,p_i\,e\otimes
e\,h(p_i)\,e)+\sum_{i=1}^n \lambda_i (e\,s(p_i)\,e\otimes e\,p_i\,e) +\\
+\sum_{i=1}^n \lambda_i \sum_{j=2}^{n_i} e(\alpha^i_{n_i}\cdots
\alpha^i_j)e\otimes e(\alpha^i_{j-1}\cdots\alpha^i_1)e.
\end{multline*}
 It follows that, for each path $q$ in $Q$, $eqe=q$ if $q$ starts and
 ends at vertices in $X$, and zero otherwise. Thus,
\begin{multline*}\Delta_{eCe}(p)=\sum_{i=1}^n
\lambda_i (p_i\otimes h(p_i))+ \sum_{i=1}^n \lambda_i (s(p_i)\otimes
p_i)+\\ +\sum_{i=1}^n \lambda_i \sum_{j=2}^{r_i} q^i_{r_i}\cdots
q^i_j\otimes q^i_{j-1}\cdots q^i_1.\end{multline*}
 Hence, this is a
linear combination of linearly independent vectors of the vector
space $eCe\otimes eCe$, so $\Delta_{eCe}(p)=y\otimes p+p\otimes x$
if and only if we have
\begin{enumerate}[$(a)$]
\item $h(p_i)=x$ for all $i=1,\ldots ,n$;
\item $s(p_i)=y$ for all $i=1,\ldots ,n$;
\item $\sum_{i=1}^n \lambda_i \sum_{j=2}^{r_i} q^i_{r_i}\cdots q^i_j\otimes
q^i_{j-1}\cdots q^i_1=0$.
\end{enumerate}

The condition $(c)$ is satisfied if and only if $r_i=1$ for all
$i=1,\ldots , n$. Therefore $\Delta_{eCe}(p)=y\otimes p+p\otimes x$
if and only if $p_i$ is a cell from $x$ to $y$ for all $i=1,\ldots
,n$. Thus the vector space of all non-trivial $(x,y)$-primitive
elements is $K\mathcal{C}ell^Q_X(x,y)\cap C$.

\begin{prop}\label{coallocaliz1} Let $C$ be an admissible subcoalgebra of the path
coalgebra $KQ$ of a quiver $Q$. Let $e_X$ be the idempotent in $C^*$
associated to a subset of vertices $X$. Then the localized coalgebra
$e_XCe_X$ is an admissible subcoalgebra of the path coalgebra
$KQ^{e_X}$, where $Q^{e_X}$ is the quiver whose set of vertices is
$(Q^{e_X})_0=X$ and the number of arrows from $x$ to $y$ is $\dim_K
K\mathcal{C}ell^Q_X(x,y)\cap C$ for all $x,y\in X$.
\end{prop}
\begin{pf} It follows from the above discussion.
\end{pf}

\begin{cor}
Let $Q$ be a quiver and $e_X$ be the idempotent of $(KQ)^*$
associated to a subset of vertices $X$. Then the localized coalgebra
$e_X(KQ)e_X$ is the path coalgebra $KQ^{e_X}$, where
$Q^{e_X}=(X,\mathcal{C}ell^Q_X)$.
\end{cor}
\begin{pf}
By \cite{gabriel}, the global dimension of $e_XCe_X$ is less or
equal than the global dimension of $C$. Then $e_XCe_X$ is
hereditary, i.e., it is the path coalgebra of its Gabriel quiver,
that is isomorphic to the quiver $Q^{e_X}$, by Proposition
\ref{coallocaliz1}.
\end{pf}

\begin{exmp}\label{example1} Let $Q$ be the quiver
\[ \xy *+{\txt{\scriptsize $x_{_1}$}}

  \PATH ~={**\dir{}} '(4,0)*+{\circ} ~={**\dir{-}} ~<{|>*\dir{>}} '(19,10)*+{\bullet}^-{\alpha_{_1}}
  ~<{|>*\dir{>}} (34,0)*+{\circ}^-{\alpha_{_2}}

  \PATH ~={**\dir{}} '(4,0)*+{\circ} ~={**\dir{-}} ~<{|>*\dir{>}}
  '(19,-10)*+{\circ}_-{\alpha_{_3}} ~<{|>*\dir{>}}
  (34,0)*+{\circ}_-{\alpha_{_4}}

  \PATH ~={**\dir{}} '(19,13)*+{\txt{\scriptsize $x_{_2}$}} '(19,-13)*+{\txt{\scriptsize
  $x_{_3}$}} (38,0)*+{\txt{\scriptsize $x_{_4}$}}

  \endxy \] and $C$ be the admissible subcoalgebra generated by
  $\alpha_2\alpha_1+\alpha_4\alpha_3$. Let us consider $X=\{x_1,
  x_3, x_4\}$. Then $e_XCe_X$ is the path coalgebra of the quiver
  $$Q^{e_X}\equiv \xy \xymatrix{\circ \ar[r]^-{\overline{\alpha_3}} & \circ
  \ar[r]^-{\overline{\alpha_4}} & \circ}
\POS (0,-3)*+{\txt{\scriptsize $x_{_1}$}}
\POS (13,-3)*+{\txt{\scriptsize $x_{_3}$}}
\POS (26,-3)*+{\txt{\scriptsize $x_{_4}$}}
\endxy
  $$
Here, the element $\alpha_2\alpha_1+\alpha_4\alpha_3$ corresponds to
the composition of the arrows $\overline{\alpha_3}$ and
$\overline{\alpha_4}$ of $Q^{e_X}$.

On the other hand, if we consider the path coalgebra $KQ$, the
quiver associated to $e_X(KQ)e_X$ is the following:
 $$Q_{e_X}\equiv \xy \xymatrix{\circ \ar[r]_-{\overline{\alpha_3}}
 \ar@/^8pt/[rr]^-{\beta\equiv \alpha_2\alpha_1}& \circ
  \ar[r]_-{\overline{\alpha_4}} & \circ}
\POS (0,-3)*+{\txt{\scriptsize $x_{_1}$}}
\POS (13,-3)*+{\txt{\scriptsize $x_{_3}$}}
\POS (26,-3)*+{\txt{\scriptsize $x_{_4}$}}
\endxy
  $$ and $\alpha_2\alpha_1+\alpha_4\alpha_3$ corresponds
to the element $\beta+\overline{\alpha_4}\hspace{0.1cm}
\overline{\alpha_3}$.
\end{exmp}

\begin{rem} As in the previous example, it is worth pointing out that if
$C$ is a proper admissible subcoalgebra of a path coalgebra $KQ$,
then we may consider two quivers: the quiver $Q^e$ defined above and
the quiver $Q_e$ such that $e(KQ)e\cong KQ_e$. Clearly, $Q^e$ is a
subquiver of $Q_e$ and the differences appear in the set of arrows.
\end{rem}

Let us now restrict our attention to the colocalizing subcategories
of $\M^C$. For the convenience of the reader we introduce the
following notation:

 We say $p=\alpha_n\cdots\alpha_2\alpha_1$ is a \emph{$h(p)$-tail} in $Q$ relative to $X$ if
$I_p\cap X=\{h(p)\}$ and $s(\alpha_i)\notin X$ for all $i=1,\ldots
,n$. If there is no confusion we simply say that $p$ is a tail.
Given a vertex $x\in X$, we denote by $\mathcal{T}ail^Q_X(x)$ the
set of all $x$-tails in $Q$ relative to $X$.

\begin{lem} Let $Q$ be a quiver and $X\subseteq Q_0$ be a subset of
vertices. Given a path $p$ in $Q$ such that $h(p)\in X$ and
$s(p)\notin X$, then $p$ has a unique decomposition $p=t q_r\cdots
q_1$, where $t, q_1,\ldots ,q_r$ are subpaths of $p$ such that $t\in
\mathcal{T}ail^Q_X(h(t))$ and $q_i\in\mathcal{C}ell^Q_X$ for all
$i=1,\ldots ,r$.
\end{lem}
\begin{pf}
It is straightforward.
\end{pf}

Let $p$ be a path in $Q$ such that $h(p)\in X\subseteq Q_0$ and
$s(p)\notin X$, we call the \emph{tail decomposition} of $p$
relative to $X$ to the decomposition given in the above lemma. We
say that $t$ is the tail of $p$ relative to $X$ if $p=t q_r\cdots
q_1$ is the tail decomposition of $p$ relative to $X$.

Assume that $C$ is an admissible subcoalgebra of a path coalgebra
$KQ$ and let $\{S_x\}_{x\in Q_0}$ be a complete set of pairwise non
isomorphic indecomposable simple right $C$-comodules. We recall that
a right $C$-comodule $M$ is \emph{quasifinite} if and only if
$\Hom_{C}(S_x,M)$ has finite dimension for all $x\in Q_0$. Let $x\in
Q_0$ and $f$ be a linear map in $\Hom_{C}(S_x,M)$. Then $\rho_M
\circ f=(f\otimes I)\circ \rho_{S_x}$, where $\rho_M$ and
$\rho_{S_x}$ are the structure maps of $M$ and $S_x$ as right
$C$-comodules, respectively. In order to describe $f$, since
$S_x=Kx$, it is enough to choose the image for $x$. Suppose that
$f(x)=m\in M$. Since $(\rho_M f)(x)=((f\otimes I)\rho_{S_x})(x)$, we
obtain that $\rho_M(m)=m\otimes x$. Therefore
$$M_x:=\Hom_{C}(S_x,M)\cong \{\text{$m\in M$ such that
$\rho_M(m)=m\otimes x$}\},$$ as $K$-vector spaces, and $M$ is
quasifinite if and only if $M_x$ has finite dimension for all $x\in
Q_0$.

Our aim is to establish when a localizing subcategory $\T_e$ is
colocalizing, or equivalently, by \cite{jmnr}, when $eC$ is a
quasifinite right $eCe$-comodule. We recall that the structure of
$eC$ as right $eCe$-comodule is given by $\rho_{eC}(p)=\sum_{(p)}
ep_{(1)}\otimes ep_{(2)} e$ if $\Delta_{KQ}(p)=\sum_{(p)}
p_{(1)}\otimes p_{(2)}$, for all $p\in eC$. It is easy to see that
$eC$ has a decomposition $eC=\bigoplus_{a\in X,b\in Q_0} C_{ab}$, as
vector space, that is, the elements of $eC$ are linear combinations
of paths which start at vertices in $X$.

\begin{prop}\label{eCquasifinite}
Let $C$ be an admissible subcoalgebra of a path coalgebra $KQ$, let
$X$ be a subset of $Q_0$ and let $e_X$ be the idempotent of $C^*$
associated to $X$. The following conditions are equivalent:

\begin{enumerate} [$(a)$]
\item The localizing subcategory $\mathcal{T}_X$ of $\mathcal{M}^C$ is
colocalizing.

\item $e_XC$ is a quasifinite right $e_XCe_X$-comodule.

\item $\dim_K K\mathcal{T}ail^Q_X(x)\cap C$ is finite for all
$x\in X$.
\end{enumerate}
\end{prop}
\begin{pf}
By the arguments mentioned above, it is enough to prove that
$(e_XC)_x=K\mathcal{T}ail^Q_X(x)\cap C$. For simplicity we write $e$
instead of $e_X$.

Let $p=\sum_{i=1}^n \lambda_i t_i\in C$ be a $K$-linear combination
of $x$-tail such that $t_i=\alpha^i_{r_i}\cdots \alpha^i_1$ ends at
$y_i$ for all $i=1,\ldots n$. Then,
$$\Delta_{KQ}(p)=p\otimes x+\sum_{i=1}^n
y_i\otimes \lambda_it_i+\sum_{i=1}^n\lambda_i \sum_{j=2}^{r_i}
\alpha^i_{r_i}\cdots \alpha^i_j\otimes \alpha^i_{j-1}\cdots
\alpha^i_1,$$ and then,
\begin{multline*}\rho_{eC}(p)=e\, p\otimes
e\, x\, e+\sum_{i=1}^n\lambda_i e\, y_i\otimes e\, t_i\, e+\\
+\sum_{i=1}^n\lambda_i\sum_{j=2}^{r_i} e(\alpha^i_{r_i}\cdots
\alpha^i_j)\otimes e(\alpha^i_{j-1}\cdots \alpha^i_1)e=p\otimes x
\end{multline*}
because $\alpha^i_j$ ends at a point not in $X$ for all $j=1,\ldots
,r_i$ and $i=1,\ldots n$. Thus $p\in (eC)_x$.

Conversely, consider an element $p=\sum_{i=1}^n \lambda_i
p_i+\sum_{k=1}^m \mu_k q_k\in (eC)_x$, where $s(p_i)\in X$ for all
$i=1,\ldots , n$, and $s(q_k)\notin X$ for all $k=1,\ldots ,m$.
Moreover, let us suppose that $p_i=\overline{p}^i_{r_i}\cdots
\overline{p}^i_1$ is the cellular decomposition of $p_i$ relative to
$X$ for all $i=1, \ldots ,n$, and $q_k=t_k
\overline{q}^k_{s_k}\cdots \overline{q}^k_1$ is the tail
decomposition of $q_k$ relative to $X$ for all $k=1,\ldots ,m$.
Then,
\begin{multline*}
\rho_{eC}(p)=\sum_{i=1}^n \lambda_i \sum_{j=2}^{r_i}
\overline{p}^i_{r_i}\cdots \overline{p}^i_j\otimes
\overline{p}^i_{j-1}\ldots \overline{p}^i_1+\sum_{i=1}^n\lambda_i
s(p_i)\otimes p_i+\sum_{i=1}^n \lambda_i p_i\otimes h(p_i)+\\
+\sum_{k=1}^m \mu_k t_k\otimes q_k+\sum_{k=1}^m \mu_k
\sum_{l=2}^{s_k} t_k\overline{q}^k_{s_k}\cdots
\overline{q}^k_l\otimes \overline{q}^k_{l-1}\ldots
\overline{q}^k_1+\sum_{k=1}^m \mu_k q_k\otimes h(q_k).
\end{multline*}
A straightforward calculation shows that if $\rho_{eC}(p)=p\otimes
x$ then $n=0$, $h(q_k)=x$ and $s_k=0$ for all $k=1,\ldots ,m$.
Therefore $p\in K\mathcal{T}ail^Q_X(x) \cap C$ and the proof is
finished.
\end{pf}

\begin{cor}\cite{jmnr}
Let $Q$ be a quiver and $X$ be a subset of $Q_0$. Then the following
conditions are equivalent:
\begin{enumerate}[$(a)$]
\item The localizing subcategory $\T_X$ of $\mathcal{M}^C$ is
colocalizing.
\item The localizing subcategory $\T_X$ of $\mathcal{M}^C$ is
perfect colocalizing.
\item $\mathcal{T}ail^Q_X(x)$ is a finite set for all $x\in X$. That
is,
 there is at most a finite number of paths starting at the
same point whose only vertex in $X$ is the first one.
\end{enumerate}
\end{cor}
\begin{pf} Since $\mathcal{T}ail^Q_X(x)$ is a basis of the vector
space $K\mathcal{T}ail^Q_X(x) \cap C$, applying Proposition
\ref{eCquasifinite}, we get that $(a) \Leftrightarrow (c)$. The
statements $(a)$ and $(b)$ are equivalent by \cite[Theorem
4.2]{blas2}.
\end{pf}

\begin{exmp} Consider the quiver $Q$
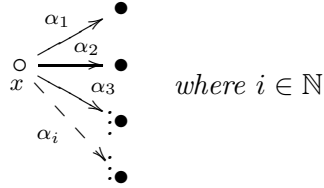
$$\xy \xymatrix@R=8pt{  &   \bullet \\
 \circ \ar[ru]^-{\alpha_1} \ar[r]^(0.66){\alpha_2}\ar[rd]^(0.7){\alpha_3} \ar@{-->}[rdd]_{\alpha_i}&   \bullet\\
       &   \bullet\\
       &   \bullet}
\POS (-0.4,-10)*+{\txt{\scriptsize $x$}} \POS (11.8,-14)*+{\vdots}
\POS (11.8,-20)*+{\vdots}
  \POS (30,-10)*+\txt{{\small where $i\in
\mathbb{N}$}}\endxy$$ and the subset $X=\{x\}$. Then
$\mathcal{T}ail^Q_X(x)=\{\alpha_i\}_{i\in \mathbb{N}}$ is an
infinite set and the localizing subcategory $\T_X$ is not
colocalizing.
\end{exmp}

\begin{rem} If $C$ is finite dimensional or, more generally, if the set $Q_0\backslash X$
 is finite and $Q$ is acyclic,
 every localizing subcategory is
colocalizing.
\end{rem}

\begin{rem} The reader should observe that all results and proofs stated in
this section remains valid for an arbitrary field $K$.
\end{rem}

\section{Comodule types of coalgebras}\label{deftame}

Let us recall from \cite{simson1} and \cite{simson2} the comodule
types of a basic (pointed) coalgebra.

 Let $C$ be a basic coalgebra such that $C_0=\oplus_{i\in I_C}
S_i$. For every finite dimensional right $C$-comodule $M$ we
consider the \emph{length vector} of $M$, $\length M=(m_i)_{i\in
I_C}\in \mathbb{Z}^{(I_C)}$, where $m_i\in \mathbb{N}$ is the number
of simple composition factors of $M$ isomorphic to $S_i$. In
\cite{simson1} it is proved that the map $M\mapsto \length M$
extends to a group isomorphism $K_0(C)\longrightarrow
\mathbb{Z}^{(I_C)}$, where $K_0(C)$ is the \emph{Grothendieck group}
of $C$.

Let $R$ be a $K$-algebra. By a $R$-$C$-bimodule we mean a $K$-vector
space $L$ endowed with a left $R$-module structure $\cdot: R\otimes
L \rightarrow L$ and a right $C$-comodule structure $\rho:L
\rightarrow L\otimes C$ such that $\rho_L(r\cdot x)=r\cdot
\rho_L(x)$. We denote by $_{R}\mathcal{M}^C$ the category of
$R$-$C$-bimodules.

\begin{lem} Let $L$ be a $R$-$C$-bimodule and $e\in C^*$ be an idempotent. Then
$eL$ is a $R$-$eCe$-bimodule and we have a functor
$$T=e(-): _{R}\!\mathcal{M}^C\longrightarrow _{R}\!\mathcal{M}^{eCe}.$$
\end{lem}
\begin{pf}
The above compatibility property yields the following equalities:
\begin{equation}\label{bicomodule} \sum_{(r\cdot
x)} (r\cdot x)_{(0)}\otimes (r\cdot x)_{(1)}=\rho_L(r\cdot x)=r\cdot
\rho_L(x)=\sum_{(x)} r \cdot x_{(0)}\otimes x_{(1)} ,
\end{equation}
for each element $r\in R$ and $x\in L$. Now, we have
$$\begin{array}{rl} r\cdot (e\cdot x)& =r\cdot (\sum_{(x)} x_{(0)}e(x_{(1)}))
\\ & =\sum_{(x)} r\cdot x_{(0)}e(x_{(1)})\\ & =(I\otimes e)(\sum_{(x)}
r\cdot x_{(0)}\otimes x_{(1)}) \\ & \overset{(1)}{=}(I\otimes e)\rho_L(r\cdot x)\\
&=e\cdot (r\cdot x).\end{array}$$ Then $eL$ has an structure of left
$R$-module and right $eCe$-comodule, and we have the compatibility
property.
$$\begin{array}{rl} r\cdot \rho_{eL}(e\cdot x)& =r\cdot (\sum_{(x)} e\cdot x_{(0)}\otimes e\cdot x_{(1)}\cdot e)
\\ &=\sum_{(x)} r\cdot (e\cdot x_{(0)})\otimes e\cdot x_{(1)}\cdot e\\
& =\sum_{(x)} e\cdot (r\cdot x_{(0)})\otimes e\cdot x_{(1)}\cdot e
\\ & \overset{(1)}{=} \sum_{(r\cdot x)} e \cdot
(r\cdot x)_{(0)}\otimes e\cdot (r\cdot x)_{(1)}\cdot e
\\ &=\rho_{eL}(e\cdot (r\cdot x)) \\& =\rho_{eL}(r\cdot (e\cdot
x)).  \end{array}$$
\end{pf}

We recall from \cite{simson1} and \cite{simson2} that a
$K$-coalgebra $C$ over an algebraically closed field $K$ is said to
be of \emph{tame comodule type} (tame for short) if for every $v\in
K_0(C)$ there exist $K[t]$-$C$-bimodules $L^{(1)},\ldots ,
L^{(r_v)}$, which are finitely generated free $K[t]$-modules, such
that all but finitely many indecomposable right $C$-comodules $M$
with $\length M=v$ are of the form $M\cong K^1_\lambda\otimes_{K[t]}
L^{(s)}$, where $s\leq r_v$, $K^1_\lambda=K[t]/(t-\lambda)$ and
$\lambda \in K$. If there is a common bound for the numbers $r_v$
for all $v\in K_0(C)$, then $C$ is called \emph{domestic}.

If $C$ is a tame coalgebra then there exists a \emph{growth
function} $\mu_C^1:K_0(C)\rightarrow \mathbb{N}$ defined as
$\mu_C^1(v)$ to be the minimal number $r_v$ of $K[t]$-$C$-bimodules
$L^{(1)},\ldots , L^{(r_v)}$ satisfying the above conditions, for
each $v\in K_0(C)$. $C$ is said to be of \emph{polynomial growth} if
there exists a formal power series
\[G(t)=\sum_{m=1}^{\infty} \sum_{j_1,\ldots ,j_m\in I_C}
g_{j_1,\ldots,j_m}t_{j_1} \ldots t_{j_m}\] with $t=(t_j)_{j\in I_C}$
and non-negative coefficients $g_{j_1,\ldots,j_m}\in \mathbb{Z}$
such that $\mu_C^1(v)\leq G(v)$ for all $v=(v(j))_{j\in I_C}\in
K_0(C)\cong \mathbb{Z}^{(I_C)}$ such that $\|v\|:=\sum_{j\in I_C}
v(j)\geq 2$. If $G(t)=\sum_{j\in I_C} g_jt_j$, where $g_j\in
\mathbb{N}$, then $C$ is called of \emph{linear growth}. If
$\mu_C^1$ is zero we say that $C$ is of discrete comodule type.

Let $Q$ be the quiver $ \xymatrix{ \circ \ar@<0.6ex>[r]
\ar@<-0.4ex>[r] \ar@<0.1ex>[r] & \circ }$ and $KQ$ be the path
algebra of the quiver $Q$. Let us denote by $\M^f_{KQ}$ the category
of finite dimensional right $KQ$-modules. A $K$-coalgebra $C$ is of
\emph{wild comodule type} (wild for short) if there exists an exact
and faithful $K$-linear functor $F: \M^f_{KQ} \rightarrow \M_f^C$
that respects isomorphism classes and carries indecomposable right
$KQ$-modules to indecomposable right $C$-comodules.

\section{Localization and tame comodule type}\label{sectiontame}

This section and the subsequent are devoted to study the relation
between the comodule type of a coalgebra and its localized
coalgebras. Let $e$ be an idempotent in the dual algebra $C^*$. We
denote by $I_e=\{ i\in I_C \hspace{0.1cm} | \hspace{0.1cm}
eS_i=S_i\}$ and $\K=\{S_i\}_{i\in I_e}$. Let us analyze the behavior
of the length vector under the action of the quotient functor.

\begin{lem}\label{longitud}
Let $C$ be a coalgebra and $e\in C^*$ be an idempotent.
\begin{enumerate}[$(a)$]
\item If $L$ is a
finite dimensional right $C$-comodule, then $(\length L)_i=(\length
eL)_i$ for all $i\in I_e$.
\item The following diagram is commutative
\[
\xymatrix@R=15pt { \M_f^C \ar[r]^-{e(-)}
\ar[dd]_-{\length} & \M^{eCe}_f \ar[dd]^-{\length}\\
   & & \\
  K_0(C) \ar[r]^-{f} & K_0(eCe)}
\]
where $f$ is the projection from $K_0(C)\cong\mathbb{Z}^{(I_C)}$
onto $K_0(eCe)\cong \mathbb{Z}^{(I_e)}$.
\end{enumerate}
\end{lem}
\begin{pf}
$(a)$ Let $\K=\{S_i\}_{i\in I_e}$. Let us consider $0\subset L_1
\subset L_2\subset\cdots \subset L_{n-1}\subset L_n$ a composition
series for $L$. Then, we obtain the inclusions $0\subseteq eL_1
\subseteq eL_2\subseteq\cdots \subseteq eL_{n-1}\subseteq eL_n$.
Since $e(-)$ is an exact functor, $eL_j/eL_{j-1}\cong
e(L_j/L_{j-1})=eS_j$, where $S_j$ is a simple $C$-comodule for all
$j=1,\ldots ,n$. But $eS_j=S_j$ if $S_j\in \K$ and zero otherwise,
see \cite{jmnr}. Thus $(\length L)_i=(\length eL)_i$ for all $i\in
I_e$. It is easy to see that $(b)$ follows from $(a)$.

\end{pf}

\begin{cor} For any finite dimensional right $C$-comodule $M$, $\len(eM)\leq \len(M)$.
\end{cor}
\begin{pf}
By Lemma \ref{longitud}, $$\len(eM)=\sum_{i\in I_e} (\length
eM)_i\leq \sum_{i\in I_C} (\length M)_i=\len(M). $$
\end{pf}

Let us now consider the converse problem, that is, take a right
$eCe$-comodule $N$ whose length vector is known, which is the length
vector of $S(N)$? Since, in general, the functor $S$ does not
preserve finite dimensional comodules (see \cite{navarro}), we have
to impose some extra conditions. We start with a simple case.

\begin{lem}\label{lemmalength} Let $N$ be a finite dimensional right $eCe$-comodule
with $\length N=(v_i)_{i\in I_e}$. Suppose that $S(S_i)=S_i$ for
all $i\in I_e$ such that $v_i\neq 0$. Then $$\left (\length S(N) \right )_i=\left\{%
\begin{array}{ll}
    v_i, & \text{if $i\in I_e$} \\
    0, & \text{if $i\in I_C\backslash I_e$} \\
\end{array}%
\right.$$
\end{lem}
\begin{pf}
Let $0\subset N_1 \subset N_2\subset\cdots \subset N_{n-1}\subset
N_n=N$ be a composition series for $N$. Since $S$ is left exact, we
have the chain of right $C$-comodules $$0\subset S(N_1) \subset
S(N_2)\subset\cdots \subset S(N_{n-1})\subset S(N_n)=S(N).$$ Now,
for each $j=0,\ldots ,n-1$, we consider the short exact sequence
$$\xymatrix{0 \ar[r]& N_j \ar[r]^-{i} & N_{j+1} \ar[r]^-{p} & S_{j+1} \ar[r]& 0}$$ and applying the
functor $S$ we have the short exact sequence
$$\xymatrix{0 \ar[r]& S(N_j) \ar[r]^-{S(i)} & S(N_{j+1}) \ar[r]^-{S(p)} & S(S_{j+1})=S_{j+1}}$$
This sequence is exact since $S(p)$ is non-zero (otherwise $S(i)$ is
bijective and then $i$ so is). Thus $S(N_{j+1})/S(N_j)\cong S_{j+1}$
and the chain is a composition series of $S(N)$.
\end{pf}

\begin{lem}
Let $C$ be a $K$-coalgebra and $R$ be a $K$-algebra. Suppose that
$N$ is a $R$-$C$-bimodule, $M$ is a right $R$-module and $f$ is an
idempotent in $C^*$. Then $f( M \otimes_R N)\simeq M\otimes_R f N$.
\end{lem}
\begin{pf}
Let us suppose that the right $C$-comodule structure of $N$ is the
map $\rho_N$. Then $M\otimes_R N$ is endowed with a structure of
right $C$-comodule given by the map $I\otimes_R \rho_N: M\otimes_R
N\rightarrow M\otimes_R N\otimes C$ defined by $$m\otimes_R n\mapsto
m\otimes_R \left ( \sum_{(n)} n_{(0)}\otimes n_{(1)} \right
)=\sum_{(n)} (m\otimes_R n_{(0)}\otimes n_{(1)}),$$ for all $m\in M$
and $n\in N$.

Therefore $f\cdot(m\otimes_R n)=\sum_{(n)} m\otimes_R n_{(0)}\otimes
f(n_{(1)})=\sum_{(n)} m\otimes_R n_{(0)}f(n_{(1)})=m\otimes_R
\sum_{(n)}n_{(0)}f(n_{(1)})=m\otimes_R f\cdot n$ for all $m\in M$
and $n\in N$. Thus $f( M \otimes_R N)\simeq M\otimes_R f N$.
\end{pf}

\begin{prop}\label{proptame} Let $v=(v_i)_{i\in I_e}\in K_0(eCe)$ such that
$S(S_i)=S_i$ for all $i\in I_e$ with $v_i\neq 0$. If $C$
satisfies the tameness condition for $\overline{v}\in K_0(C)$ given by $$(\overline{v})_i=\left\{%
\begin{array}{ll}
    v_i, & \text{if $i\in I_e$} \\
    0, & \text{if $i\in I_C\backslash I_e$} \\
\end{array}%
\right.$$    then $eCe$ satisfies the tameness condition for $v$,
and $\mu_{eCe}^1(v)\leq \mu_C^1(\overline{v})$.
\end{prop}
\begin{pf} By hypothesis, there exist $K[t]$-$C$-bimodules
$L^{(1)}, L^{(2)}, \ldots , L^{(r_{\overline{v}})}$, which are
finitely generated free $K[t]$-modules, such that all but finitely
many indecomposable right $C$-comodules $M$ with $\length
M=\overline{v}$ are of the form $M\cong K^1_\lambda\otimes_{K[t]}
L^{(s)}$, where $s\leq r_{\overline{v}}$,
$K^1_\lambda=K[t]/(t-\lambda )$ and $\lambda \in K$. Consider the
$K[t]$-$eCe$-bimodules $eL^{(1)},\ldots , eL^{(r_{\overline{v}})}$.
Obviously, they are finitely generated free as left $K[t]$-modules.
Let now $N$ be a right $eCe$ comodule with $\length N=v$. By Lemma
\ref{lemmalength}, $\length S(N)=\overline{v}$ and therefore $S(N)
\cong K^1_\lambda\otimes_{K[t]} L^{(s)}$ for some $s\leq
r_{\overline{v}}$ and some $\lambda \in K$ (since $S$ is an
embedding, there are only finitely many $eCe$-comodules $N$ such
that $S(N)$ is not of the above form). Then, by the previous lemma,
$N\cong eS(N)\cong  K^1_\lambda\otimes_{K[t]} eL^{(s)}$. Thus $eCe$
satisfies the tameness condition for $v$ and, furthermore,
$\mu_{eCe}^1(v)\leq \mu_C^1(\overline{v})$.
\end{pf}

Following \cite{jmnr}, the idempotent $e\in C^*$ is said to be
\emph{left (right) semicentral} if $eCe=eC$ ($eCe$=$Ce$).

\begin{cor}
Let $C$ be a coalgebra and $e\in C^*$ be a right semicentral
idempotent. If $C$ is tame (of polynomial growth, of linear growth,
domestic, discrete) then $eCe$ is tame (of polynomial growth, of
linear growth, domestic, discrete).
\end{cor}
\begin{pf} By \cite{navarro}, $e$ is right semicentral if and only of $S(S_i)=S_i$ for all
$i\in I_e$. Then the statement follows from the former result.
\end{pf}

The underlying idea of the proof of Proposition \ref{proptame} goes
as follows: if we control the $C$-comodules whose length vector is
obtained as $\length S(T)$ for some $eCe$-comodule $T$ such that
$\length T=v$ (under the assumption of Proposition \ref{proptame}
the unique length vector obtained in such a way is exactly
$\overline{v}$), then we may control the $eCe$-comodules of length
$v$. In the following result we generalize Proposition
\ref{proptame} applying the aforementioned method.

For the convenience of the reader we introduce the following
notation. To any vector $v\in K_0(eCe)\cong \mathbb{Z}^{(I_e)}$ we
associate the set $\Omega_v=\{v^{(\beta)}\}_{\beta\in B}$ of all
vectors in $K_0(C)\cong \mathbb{Z}^{(I_C)}$ such that each
$v^{(\beta)}=\length S(N)$ for some $eCe$-comodule $N$ such that
$\length N=v$.
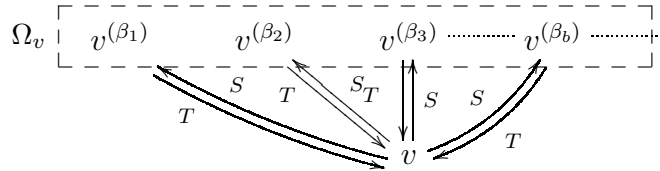
$$\xy\xymatrix@R=30pt{  v^{(\beta_1)} \ar@/_5pt/@<-1.2ex>[rrd]_(0.3){\text{\tiny $T$}}  &    v^{(\beta_2)}
\ar@<-0.6ex>[rd]
& v^{(\beta_3)}\ar@{.}[r]\ar@<-0.4ex>[d]
& v^{(\beta_b)} \ar@{.}[r] \ar@/^8pt/@<0.7ex>[ld]^-{\text{\tiny $T$}} & \\
  &    &  v  
  \ar@/^5pt/@<0.5ex>[llu]_(0.6){\text{\tiny $S$}}
   \ar[lu]_-{\text{\tiny $S$}} \ar@<-0.3ex>[u]_-{\text{\tiny $S$}}
   \ar@/_8pt/[ru]^-{\text{\tiny $S$}} &
  }
\POS (-8,-4) \ar@{--} (70,-4)

\POS (-8,-4) \ar@{--} (-8,4)

\POS (-8,4) \ar@{--} (70,4)

\POS (70,-4) \ar@{--} (70,4)

\POS (-12,0)*+{\txt{\small $\Omega_v$}}

\POS (33,-7.8)*+{\text{\tiny $T$}}

\POS (22,-7.8)*+{\text{\tiny $T$}}
\endxy
  $$

\begin{prop}\label{tamegeneralized} Let $v\in K_0(eCe)$ such that $\Omega_v$ is
finite. If $C$ satisfies the tameness condition for each vector
$v^{(\beta)}\in \Omega_v$ then $eCe$ satisfies the tameness
condition for $v$.
\end{prop}
\begin{pf}
Consider the set of all $K[t]$-$C$-bimodules associated to all
vectors $v^{(\beta)}\in \Omega_v$, namely
$\mathcal{L}=\{L^{(j)}_{\beta}\}_{\beta\in \Omega_v,j=1, \ldots
,r_\beta}$. By hypothesis, this is a finite set and then
$T(\mathcal{L})$ so is. We proceed analogously to Proposition
\ref{proptame} and the result follows.
\end{pf}

Given two vectors $v=(v_i)_{i\in I_C}, w=(w_i)_{i\in I_C} \in
K_0(C)$, we say that $v\leq w$ if $v_i\leq w_i$ for all $i\in I_C$.

\begin{lem}\label{omegafinito} Let $v=(v_i)_{i\in I_e}\in K_0(eCe)$ such that
 $S(S_i)$ is a finite dimensional right $C$-comodule for all
 $i\in I_e$ with $v_i\neq 0$. Then
 $\Omega_v$ is a finite set.
 \end{lem}
 \begin{pf} Let $N$ be a right $eCe$-comodule such that $\length N=v$.
 Consider a composition series for $N$, $0\subset N_1 \subset
N_2\subset\cdots \subset N_{n-1}\subset N_n=N$. Since $S$ is left
exact, we have the chain of right $C$-comodules $$0\subset S(N_1)
\subset S(N_2)\subset\cdots \subset S(N_{n-1})\subset S(N_n)=S(N).$$

Then, for each $j=1,\ldots ,n$, we have an exact sequence
$$\xymatrix{0 \ar[r] &S(N_{j-1})\ar[r]&  S(N_j) \ar[r] & S(S_j)},$$
where $S_j$ is a simple $eCe$-comodule.

By \cite{navarro}, $\Soc(S(S_j))=S_j$ and hence $S(S_j)$ has a
composition series
$$0\subset S_j\subset S(S_j)_2 \subset \cdots\subset S(S_j)_{r-1}
\subset S(S_j).$$ Then we can complete the following commutative
diagram taking the pullback $P_i$ for $i=1,\ldots ,r-1$.

$$\xy \xymatrix@C=30pt{ S(N_j)\ar[r]\ar[r]|(0.3){|}\ar@2@{=}[d] &
  P_1\ar[d]\ar[d]|(0.27){-} \ar[r]\ar[r]|(0.8){|} & S_j\ar[d]\ar[d]|(0.3){-}
\\
S(N_j)\ar[r]\ar[r]|(0.3){|}\ar@2@{=}[d] &  P_2
\ar[r]\ar[d]\ar[d]|(0.27){-} & S(S_j)_2 \ar[r]\ar[d]\ar[d]|(0.3){-}
&
\Coker \\
\vdots \ar@2@{=}[d]& \vdots \ar[d]\ar[d]|(0.35){-} &\vdots
\ar[d]\ar[d]|(0.35){-}&\vdots
\\
S(N_j)\ar[r]\ar[r]|(0.3){|}\ar@2@{=}[d] &  P_{r-1}
\ar[r]\ar[d]\ar[d]|(0.3){-} & S(S_j)_{r-1}
\ar[r]\ar[d]\ar[d]|(0.3){-} & \Coker\\
S(N_j)\ar[r]\ar[r]|(0.3){|}&  S(N_{j+1}) \ar[r]  &  S(S_j) \ar[r] &
\Coker }
%
\endxy
$$
Consider two consecutive rows and their quotient sequence
$$\xy \xymatrix@R=15pt@C=20pt{ S(N_j)\ar[r]\ar[r]|(0.38){|}\ar@2@{=}[dd] &
  P_t\ar[dd]^-{i}\ar[dd]|(0.2){-} \ar[rr]^-{g_t}\ar[rd] &
  &
  S(S_j)_t\ar[dd]^-{i}\ar[dd]|(0.2){-}\ar[r]&
  \Coker
\\
 &                 &        \Ima g_t \ar[ru] \ar[dd]^(0.4){i}\ar[dd]|(0.2){-}  &     &         \\
S(N_j)\ar[r]\ar[r]|(0.38){|}\ar[dd] &
  P_{t+1} \ar[rr]^(0.3){g_{t+1}}\ar[dd]^-{p}\ar[dd]|(0.7){-}\ar[rd]&
 &
 S(S_j)_{t+1}\ar[r]\ar[dd]^-{p}\ar[dd]|(0.73){-} &
\Coker \\
&                 &        \Ima g_{t+1} \ar[ru] \ar[dd]^(0.4){p}\ar[dd]|(0.73){-}  &     &         \\
0 \ar[r]&
 P_{t+1}/P_t \ar[rr]^(0.3){\overline{g}} \ar[rd]&
  &
 S_k  & \\
&                 &       \Ima g_t/\Ima g_{t+1} \ar[ru]   & & }
 \POS

 \POS (13.3,-3)+*{\xymatrix@R=1pt@C=20pt{ \ar[rd]& \\ & 0}}

 \POS (13.6,-8.5)+*{\xymatrix@R=1pt@C=20pt{ \ar[rd]& \\ & 0}}

\POS (16.5,-13.5)+*{\xymatrix@C=15pt{\cong\Ima \overline{g} \ar[r]&
0 }}

\POS (47.5,-21.5)+*{\cong S(S_j)_{t+1}/S(S_j)_t}
\endxy $$
Suppose that $P_{t+1}\neq P_t$, then $P_{t+1}/ P_t\cong \Ima
\overline{g}\hookrightarrow S_k$, and thus $P_{t+1}/ P_t\cong S_k$.

Hence we have obtained a chain
$$0\subset P^1_1\subseteq \cdots \subseteq P^1_{r_1}=S(N_1)\subseteq
\cdots\subseteq S(N_{n-1})\subseteq P^n_1\subseteq\cdots \subseteq
P^n_{r_n}=S(N),$$ where the quotient of two consecutive comodules is
zero or a simple comodule.

Therefore $\length S(N)\leq \sum_{j=1}^n \length S(S_j)$ for any
right $eCe$-comodule $N$ whose $\length N=v$. Thus $\Omega_v$ is a
finite set.
\end{pf}

As a consequence of the former results we may state the following
theorem:

 \begin{thm}\label{corolariotame} Let $C$ be a coalgebra and $e\in C^*$ be an idempotent such that
 $S(S_i)$ is a finite dimensional
 right $C$-comodule for all $i\in I_e$. If $C$ is of tame (discrete) comodule type then $eCe$ is
 of tame (discrete) comodule type.
 \end{thm}
 \begin{pf} By Lemma \ref{omegafinito}, $\Omega_v$ is a finite
 set for each $v\in K_0(eCe)$. Then, by Proposition
 \ref{tamegeneralized}, $eCe$ satisfies the tameness condition for
 $v$.
 \end{pf}

\begin{rem}\label{Spreserfindim} The reader should observe that the proof
of Lemma \ref{omegafinito} also shows that the section functor $S$
preserves finite dimensional comodules if and only if $S(S_i)$ is
finite dimensional for all $i\in I_e$.
\end{rem}

In particular, the conditions of Theorem \ref{corolariotame} are
satisfied for any idempotent if $C$ is left semiperfect. A coalgebra
is said to be \emph{left semiperfect} if every finite dimensional
left comodule has a finite-dimensional projective cover, or
equivalently, if any indecomposable injective right comodule is
finite dimensional.

\begin{cor} Let $C$ be a left semiperfect coalgebra
and $e\in C^*$ be an idempotent. If $C$ is of tame (discrete)
comodule type then $eCe$ is of tame (discrete) comodule type.
\end{cor}
\begin{pf} By \cite[Proposition 3.2 and Corollary 3.3]{navarro}, for
any $i\in I_e$, $S(S_i)$ is a subcomodule of an indecomposable
injective right $C$-comodule. Therefore, if $C$ is left semiperfect,
$S(S_i)$ is finite dimensional for any $i\in I_e$. Thus the result
follows from Theorem \ref{corolariotame}.
\end{pf}

The following problem remains still open.

\begin{prob}\label{problem1}
Assume that $C$ is a coalgebra of tame comodule type and $e$ is an
idempotent in $C^*$. Is the coalgebra $eCe$ of tame comodule type?
\end{prob}
It would be interesting to know if the localization process
preserves polynomial growth, linear growth, discrete comodule type
or domesticity. It is clear that the converse result is not true as
the following example shows.
\begin{exmp}
Let us consider the quiver $$\xymatrix@C=15pt@R=10pt{ & & & \circ \ar[d] & & & &\\
 Q: & \circ \ar[r] & \circ \ar[r] &\circ \ar[r]
 &\circ \ar[r] & \circ \ar[r] &\circ \ar[r] &\circ \ar[r]
 &\circ \ar[r] & \circ
 }$$ Since its underlying graph is neither a Dynkin diagram nor an Euclidean
 graph, $KQ$ is wild, see \cite[Theorem 9.4]{simson2}. But it is easy to see that $eCe$ is of
 tame comodule type for each non-trivial idempotent
 $e\in C^*$.
\end{exmp}

\section{Split idempotents}\label{sectionwild}

Let us study the wildness of a coalgebra and its localized
coalgebras. Directly from the definition we may prove the following
proposition.

\begin{prop}\label{perfectcolocalizwild}
Let $C$ be a coalgebra and $e\in C^*$ be an idempotent which defines
a perfect colocalization. If $eCe$ is wild then $C$ is wild.
\end{prop}
\begin{pf}By hypothesis, there is an exact and faithful functor
$F:\M^f_{KQ}\rightarrow \M_f^{eCe}$, where $Q$ is the quiver $
\xymatrix{ \circ \ar@<0.6ex>[r] \ar@<-0.4ex>[r] \ar@<0.1ex>[r] &
\circ }$, which respects isomorphism classes and preserves
indecomposables. In \cite{navarro} it is proved that the
colocalizing functor $H$ restricts to a functor
$H:\M_f^{eCe}\rightarrow \M_f^C$ that preserves indecomposables and
respects isomorphism classes. $H$ is also exact by hypothesis.
Therefore the composition $HF:\M^f_{KQ}\rightarrow \M_f^{C}$ is an
exact and faithful functor that preserves indecomposables and
respects isomorphism classes. Thus $C$ is wild.

\end{pf}

A similar result may be obtained using the section functor if $S$
preserves finite dimensional comodules. For example, if $C$ is left
semiperfect.

\begin{prop}Let $C$ be a left semiperfect coalgebra and
$e\in C^*$ be an idempotent which defines a perfect localization. If
$eCe$ is wild then $C$ is wild.
\end{prop}
\begin{pf} It is analogous to the proof of Proposition
\ref{perfectcolocalizwild}.
\end{pf}

Let us now consider the following question: when is the coalgebra
$eCe$ a subcoalgebra of $C$? This is interesting for us because, by
\cite[Theorem 5.4]{simson2}, in such a case, we obtain the following
result:

\begin{prop}\label{subcoalwild} Let $C$ be a coalgebra and $e\in C^*$ be an
idempotent such that $eCe$ is a subcoalgebra of $C$. If $eCe$ is
wild then $C$ is wild.
\end{prop}

In general, we always have the inclusion $eCe\subseteq{C}$,
nevertheless the coalgebra structures may be different. This is not
the case if, for instance, $e$ is a left semicentral idempotent. In
that case, by \cite{jmnr}, $eC=eCe$ is a subcoalgebra of $C$. The
same result holds if $e$ is a right semicentral or a central
idempotent.

An idempotent $e\in{C^*}$ is said to be \emph{split} if in the
decomposition $C^*=eC^*e\oplus eC^*f\oplus{fC^*e}\oplus{fC^*f}$,
where $e+f=1$, the direct summand
$H_e:=eC^*f\oplus{fC^*e}\oplus{fC^*f}$ is a twosided ideal of $C^*$.
These elements were used by Lam in \cite{Lam}. The main result
there, see \cite[Theorem 4.5]{Lam}, assures that the following
statements are equivalent:
 \begin{enumerate}[$(a)$]
 \item $H_e$ is a twosided ideal of $C^*$.
 \item $e(C^*fC^*)e=0$.
 \item $exeye=exye$ for any $x$, $y\in{C^*}$.
  \end{enumerate}
  As a consequence, every left or right semicentral idempotent
   in $C^*$ is split. Let us characterize when
$eCe$ is a subcoalgebra of $C$.

\begin{thm}
Let $e$ be an idempotent in $C^*$. Then the following statements are
equivalent.
 \begin{enumerate}[$(a)$]
 \item $e$ is a split idempotent in $C^*$.
 \item $eCe$ is a subcoalgebra of $C$.
 \end{enumerate}
\end{thm}
\begin{pf}
Let us denote $f=1-e$. For any subspace $V\subseteq{C}$, $V$ is a
subcoalgebra of $C$ if and only if $V^{\perp}$ is a twosided ideal
of $C^*$. Then we proceed as follows in order to compute the
orthogonal of $eCe$.
\[
\begin{array}{rl}
 (eCe)^{\perp}=
 &(eC\cap{Ce})^{\perp} \\
 =&(eC)^{\perp}+(Ce)^{\perp}\\
 =&C^*f+fC^*\\
 = &eC^*f+fC^*f+fC^*e+fC^*f \\
 =&eC^*f+fC^*e+fC^*f \\
 =&H_e
\end{array}
\]
Thus $eCe$ is a subcoalgebra of $C$ if and only if $H_e$ is an ideal
of $C^*$ if and only if $e$ is a split idempotent in $C^*$.
\end{pf}

Let us give a description of the split idempotents. Suppose that $C$
is a pointed coalgebra, that is, $C$ is an admissible subcoalgebra
of a path coalgebra. We recall from \cite{jmnr} or \cite{navarro}
that left (right) semicentral idempotents can be described as
follows.

\begin{prop}
Let $C$ be an admissible subcoalgebra of a path coalgebra $KQ$ and
$e_X$ be the idempotent in $(KQ)^*$ associated to a subset
$X\subseteq Q_0$. Then:
\begin{enumerate}[$(a)$]
\item $e_X$ is left semicentral if and only if there is no arrow $y\rightarrow x$ in
$Q$ such that $y\notin X$ and $x\in X$.
\item $e_X$ is right semicentral if and only if there is no arrow $x\rightarrow y$ in
$Q$ such that $y\notin X$ and $x\in X$.
\end{enumerate}
\end{prop}

We want to give a geometric description of the split idempotents in
a similar way. In order to do this, we start giving an approach by
means of path coalgebras.

\begin{lem}\label{splitidemp}
Let $Q$ be a quiver and $e_X\in (KQ)^*$ be the idempotent associated
to a subset of vertices $X$. Then $e_X$ is split in $(KQ)^*$ if and
only if $I_p\subseteq X$ for any path $p$ in $e_X(KQ)e_X$, i.e.,
there is no cell in $Q$ relative to $X$ of length greater than one.
\end{lem}
\begin{pf}
Note that $e_X(KQ)e_X$ is a subcoalgebra of $KQ$ if and only if
$\Delta(p)\in e_X(KQ)e_X\otimes e_X(KQ)e_X$ for any path $p$ in
$e_X(KQ)e_X$.

Let $p=\alpha_n\cdots\alpha_1\in e_X(KQ)e_X$, $\Delta(p)\in
e_X(KQ)e_X\otimes e_X(KQ)e_X$ if and only if $$\sum_{j=2}^n
\alpha_n\cdots\alpha_j\otimes \alpha_{j-1}\cdots\alpha_1\in
e_X(KQ)e_X\otimes e_X(KQ)e_X.$$ Since all summands are linearly
independent, this happens if and only if $s(\alpha_i)\in X$ for all
$i=1,\ldots n-1$. That is, if and only if $I_p\subseteq X$.
\end{pf}

Note that Lemma \ref{splitidemp} asserts that the vertices
associated to a split idempotent is a convex set of vertices in the
quiver $Q$. Then it is really easy to decide whether or not
$e_X(KQ)e_X$ is a subcoalgebra of $KQ$.

\begin{exmp} Let $Q$ be the following quiver:
$$
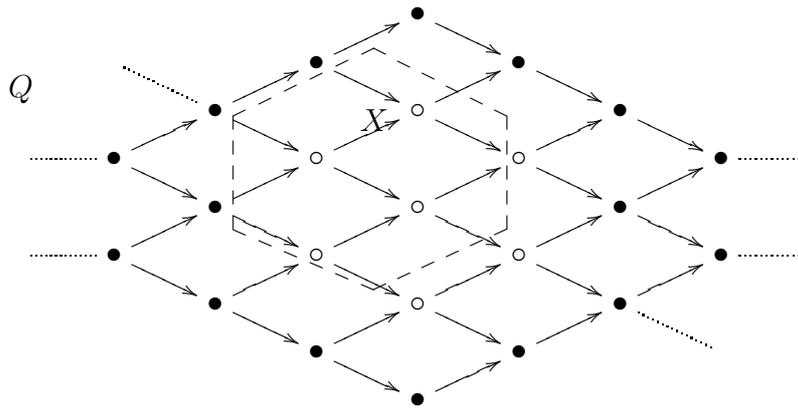
\xy
\xymatrix@R=5pt{    &       &      &   &  \bullet \ar[rd] &  &    &     &\\
    &       &      &  \bullet  \ar[rd] \ar[ru] &    &  \bullet \ar[rd]  &    &     &\\
    &       &  \bullet \ar[ru] \ar[rd] \ar@{.}[lu]   &      & \circ \ar[ru] \ar[rd]   &     & \bullet \ar[rd]   &     &\\
  \ar@{.}[r]  &  \bullet \ar[ru] \ar[rd]     &      & \circ \ar[rd] \ar[ru]     &    & \circ  \ar[rd] \ar[ru]   &    &  \bullet \ar@{.}[r]   &\\
    &       &  \bullet \ar[rd] \ar[ru]    &      & \circ \ar[rd] \ar[ru]   &     &  \bullet\ar[rd]\ar[ru]  &     &\\
  \ar@{.}[r]  &  \bullet \ar[ru] \ar[rd]     &      & \circ \ar[rd] \ar[ru]     &    &  \circ \ar[rd] \ar[ru]   &    &  \bullet\ar@{.}[r]   &\\
    &       &  \bullet \ar[ru] \ar[rd]    &      & \circ \ar[ru] \ar[rd]   &     & \bullet \ar[ru]  \ar@{.}[rd] &     &\\
  &       &      &  \bullet  \ar[ru]\ar[rd] &    &  \bullet \ar[ru]  &    &
  & \\
 &       &      &   &   \bullet \ar[ru] &   &    &&   }

\POS (64,-13.5) \ar@{--}  (64,-28.5)

\POS (28,-13.5) \ar@{--}  (28,-28.5)

\POS (64,-13.5) \ar@{--}  (46.5,-4.5)

\POS (28,-13.5) \ar@{--}  (46.5,-4.5)

\POS (64,-28.5) \ar@{--}  (46.5,-36.5)

\POS (28,-28.5) \ar@{--}  (46.5,-36.5)

\POS (46.5,-14)*+{X}

\POS (0,-10)*+{Q}
\endxy
  $$
Then the idempotent associated to the set of white vertices $X$ is a
split idempotent.
\end{exmp}

The proof of Lemma \ref{splitidemp} easily extends to pointed
coalgebras. To do that, we denote by $\mathcal{Q}=Q_0\cup Q_1\cup
\cdots \cup Q_n\cup \cdots$ the set of all paths in $Q$. Let $a$ be
an element of $KQ$. Then we can write $a=\sum_{p\in \mathcal{Q}} a_p
p$, for some $a_p\in K$. We define \emph{the path support} of $a$ to
be $\Sup(a)=\{\text{$p\in \mathcal{Q}$ $|$ $a_p\neq 0$}\}$.
Moreover, for any set $R\subseteq KQ$, we define
$\Sup(R)=\bigcup_{a\in R} \Sup(a)$.
\begin{lem}
Let $Q$ be a quiver and $C$ be an admissible subcoalgebra of $KQ$.
Let $e_X\in C^*$ be the idempotent associated to a subset of
vertices $X$. Then $e_X$ is split in $C^*$ if and only if
$I_p\subseteq X$ for any path $p$ in $\Sup(e_XCe_X)$.
\end{lem}

\begin{exmp} Let $Q$ be the quiver
$$\xy \xymatrix{ \circ \ar[r]^-{\alpha} & \circ \ar[r]^-{\beta} &
\circ} \POS (0,3)*+{\txt{\scriptsize 1}} \POS
(12,3)*+{\txt{\scriptsize 2}} \POS (25,3)*+{\txt{\scriptsize 3}}
\endxy $$
and $C$ be the admissible subcoalgebra of $KQ$ generated by $\{1, 2,
3, \alpha , \beta\}$. Then $e\equiv \{1,3\}$ is a split idempotent
because $eCe=S_1\oplus S_3$ is a subcoalgebra of $C$.
\end{exmp}

Let us finish the section with an open problem for further
development of representation theory of coalgebras.

\begin{prob}\label{problem2}
Let $C$ be a coalgebra and $e\in C^*$ be an idempotent. If $eCe$ is
of wild comodule type then $C$ is of wild comodule type.
\end{prob}

Obviously, Problem \ref{problem1} and Problem \ref{problem2} are
equivalent if the \emph{tame-wild dichotomy} for coalgebras,
conjectured by Simson in \cite{simson2}, is true.

\section{A Gabriel's theorem for coalgebras}\label{sectiongabriel}

Following \cite{simson2}, the authors intent to prove in \cite{jmn}
that every basic coalgebra, over an algebraically closed field $K$,
is the path coalgebra $C(Q,\Omega)$ of a quiver $Q$ with a set of
relations $\Omega$. That is, an analogue for coalgebras of Gabriel's
theorem. The result is proven in \cite{simson2} for the family of
coalgebras $C$ such that the Gabriel quiver $Q_C$ of $C$ is
intervally finite. Unfortunately, that proof does not hold for
arbitrary coalgebras, as a class of counterexamples given in
\cite{jmn} shows.

 Nevertheless, it is worth noting that these counterexamples are of wild comodule type
in the sense described in Section \ref{deftame}. Since, according to
Theorem \ref{criterion}, the coalgebras $C$ which are not the path
coalgebra of quiver with relations are close to be wild, we may
reformulate the problem stated in \cite[Section 8]{simson1} as
follows: \textit{any basic tame $K$-coalgebra $C$, over an
algebraically closed field $K$, is isomorphic to the path coalgebra
$C(Q,\Omega)$ of a quiver with relations $(Q,\Omega)$}. This section
is devoted to solve this problem when the Gabriel quiver $Q_C$ of
$C$ is assumed to be acyclic.

\begin{exmp} Let $Q$ be the quiver
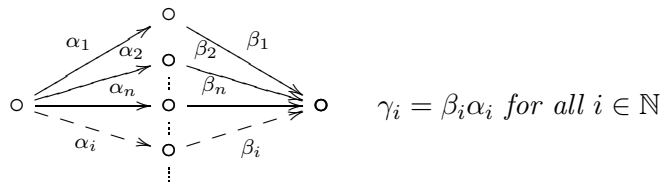
\begin{equation*} \xy *+{\circ} \PATH ~={**\dir{-}} ~<{|>*\dir{>}}
'(20,12)*+{\circ}^-{\alpha_1} ~<{|>*\dir{>}}
'(40,0)*+{\circ}^-{\beta_1} \PATH ~={**\dir{-}} ~<{|>*\dir{>}}
'(20,6)*+{\circ}^(0.8){\alpha_2} ~<{|>*\dir{>}}
'(40,0)*+{\circ}^(0.2){\beta_2} \PATH ~={**\dir{-}} ~<{|>*\dir{>}}
'(20,0)*+{\circ}^(0.7){\alpha_n} ~<{|>*\dir{>}}
'(40,0)*+{\circ}^(0.3){\beta_n} \PATH ~={**\dir{--}} ~<{|>*\dir{>}}
'(20,-6)*+{\circ}_-{\alpha_{i}} ~<{|>*\dir{>}}
'(40,0)*+{\circ}_-{\beta_{i}} \PATH  ~={**\dir{}} '(20,6)*+{\circ}
~={**\dir{.}} '(20,0)*+{\circ} '(20,-6)*+{\circ} (20,-10)
\endxy \hspace{0.5cm}\text{{\small $\gamma_i=\beta_i \alpha_i$ for all $i\in
\mathbb{N}$}}\end{equation*}
 and suppose that $H$ is the admissible subcoalgebra of $KQ$ generated by the set \mbox{$\Sigma=\{\gamma_i
 -\gamma_{i+1}\}_{i\in \mathbb{N}}$.} It is proved in \cite[Example
4.11]{jmn} that $H$ is not the path coalgebra of a quiver with
 relations. Nevertheless, $H$ is of wild comodule type because it contains
 the path coalgebra of the quiver $$ \xymatrix@R=5pt{ & \circ&   &    \circ \\
 \Gamma\equiv& &  \circ \ar[ru]^-{\alpha_1} \ar[r]^(0.7){\alpha_2} \ar[rd]_-{\alpha_3}
    \ar[ld]^-{\alpha_4} \ar[lu]_-{\alpha_5} &    \circ \\
& \circ  &   &    \circ }$$ Since $K\Gamma$ is a finite dimensional
coalgebra then we have an algebra isomorphism $(K\Gamma)^*\cong A$,
where $A$ is the path algebra of the quiver $\Gamma$, and there
exists an equivalence between the category of finite dimensional
right $A$-modules and the category of finite dimensional right
$K\Gamma$-comodules. But it is well known that $A$ is a wild algebra
and hence $K\Gamma$ is a wild coalgebra. By \cite[Theorem
5.4]{simson2}, this proves that $H$ is wild.
\end{exmp}

\begin{thm}\label{thm}
Let $K$ be an algebraically closed field, $Q$ be an acyclic quiver
and $C$ be an admissible subcoalgebra of $KQ$ which is not the path
coalgebra of a quiver with relations. Then $C$ is of wild comodule
type.
\end{thm}
\begin{pf}
By Theorem \ref{criterion}, since $C$ is not the path coalgebra of a
quiver with relations, there exists an infinite number of paths
$\{\gamma_{i}\}_{n\in \mathbb{N}}$ in $Q$ between two vertices $x$
and $y$ such that:

$\bullet$ None of them is in $C$.

$\bullet$ $C$ contains a set $\Sigma=\{\Sigma_n\}_{n\in \mathbb{N}}$
such that $\Sigma_n=\gamma_n+\sum_{j>n} a^n_j \gamma_j$, where
$a_j^n\in K$ for all $j,n\in \mathbb{N}$.
$$\xy *+{\circ} \ar @/^30pt/ ^-{\gamma_1}
(40,0)*+{\circ}
*+{\circ} \ar @/^17pt/ ^-{\gamma_2}  (40,0)*+{\circ}
*+{\circ} \ar ^(0.4){\gamma_n}  (40,0)*+{\circ}
*+{\circ} \ar @/_17pt/ _(0.6){\gamma_i}  (40,0)*+{\circ}
\POS (20,4) \ar @{.} (20,1)
\POS (20,-1) \ar @{.} (20,-5)
\POS (20,-6) \ar @{.} (20,-8)

\POS (0,-3)*+{\txt{\scriptsize $x$}}

\POS (40,-3)*+{\txt{\scriptsize $y$}}
\endxy $$
Consider $\Sup(\Sigma_1\cup \Sigma_2\cup \Sigma_3)=\{\gamma_1,
\gamma_2, \ldots , \gamma_t\}$ and $\Gamma$ the finite subquiver of
$Q$ formed by the paths $\gamma_i$ for $i=1,\ldots ,t$.

Then $D=K\Gamma\cap C$ is a finite dimensional subcoalgebra of $C$
(and an admissible subcoalgebra of $K\Gamma$) which contains the
elements $\Sigma_1$, $\Sigma_2$ and $\Sigma_3$. It is enough to
prove that $D$ is wild.

Consider the idempotent $e\in D^*$ such that $e(x)=e(y)=1$ and zero
otherwise, i.e., its associated subset of vertices is $X=\{x,y\}$.
Then, by Proposition \ref{coallocaliz1}, each $\Sigma_i$ corresponds
to an arrow from $x$ to $y$ in the quiver $\Gamma^e$, that is,
$\Gamma^e$ contains the subquiver $ \xymatrix{ \circ \ar@<0.6ex>[r]
\ar@<-0.4ex>[r] \ar@<0.1ex>[r] & \circ }$ and then $\dim_K
\Ext^1_{eDe}(S_x,S_y)\geq 3$. Thus $eDe=K\Gamma^e$ is wild by
\cite[Corollary 5.5]{simson2}. Note also that the quiver $\Gamma^e$
is of the form
\[ \xy \xymatrix@C=50pt{ \circ \ar@<1.5ex>[r]^-{\alpha_1}
\ar@<0.75ex>[r] \ar@<-1ex>[r]_-{\alpha_n} & \circ}

\POS (11,0.3)*+{\cdot} \POS (11,-0.5)*+{\cdot} \POS
(11,-1.3)*+{\cdot} \POS (0,-2.5)*+{\txt{\scriptsize $x$}} \POS
(22,-2.5)*+{\txt{\scriptsize $y$}}
\endxy
\] an then the simple right $eDe$-comodule $S_x$ is injective.

Let us prove that the localizing subcategory $\T_e\subseteq \M^D$ is
perfect colocalizing.

Since $\Gamma$ is finite and acyclic then $\dim_K
K\mathcal{T}ail^\Gamma_X(x)$ is finite and $\dim_K
K\mathcal{T}ail^\Gamma_X(y)=0$ so, by Proposition
\ref{eCquasifinite}, the subcategory $\T_e$ is colocalizing.

Let now $g$ be an element in $eC(1-e)$. Then $g$ is a linear
combination of tails starting at $x$ and then
$\rho_{eC(1-e)}(g)=g\otimes x$ (see the proof of Proposition
\ref{eCquasifinite}). Therefore $<g>\cong S_x$ as right
$eCe$-comodules. Suppose that $m=\dim_K eC(1-e)$. Hence
$eC=eCe\oplus eC(1-e)=eCe\oplus S_x^m$ and $eC$ is an injective
right $eCe$-comododule. Thus the colocalization is perfect and, by
Proposition \ref{perfectcolocalizwild}, $D$ is wild.
\end{pf}

Now we are able to extend \cite[Theorem 3.14 (c)]{simson2}, from the
case $Q$ is an intervally finite quiver and $C'\subseteq KQ$ is an
arbitrary admissible subcoalgebra to the case $Q$ is acyclic and
$C'\subseteq KQ$ is tame, as follows.

\begin{cor}[Acyclic Gabriel's theorem for coalgebras]\label{maincor}
Let $Q$ be an acyclic quiver and let $K$ be an algebraically closed
field.
\begin{enumerate}[$(a)$]
\item Any tame admissible subcoalgebra $C'$ of the path coalgebra
$KQ$ is isomorphic to the path coalgebra $C(Q,\Omega)$ of a quiver
with relation $(Q,\Omega)$.
\item The map $\Omega \mapsto C(Q,\Omega)$ defines a one-to-one
correspondence between the set of relation ideals $\Omega$ of the
path $K$-algebra $KQ$ and the set of admissible subcoalgebras $H$ of
the path coalgebra $KQ$. The inverse map is given by $H\mapsto
H^{\perp}$.
\end{enumerate}
\end{cor}\label{mainth}
\begin{pf}
Apply Theorem \ref{thm} and the weak tame-wild dichotomy proved by
Simson in \cite{simson2}.
\end{pf}

\begin{rem} It follows from the proof of Theorem \ref{thm} that if
$Q$ is acyclic, a basis of a tame admissible subcoalgebra $C$ cannot
contain three linearly independent combinations of paths with common
source and common sink and then $C=KQ_0\bigoplus \bigoplus_{x\neq y
\in Q_0} C_{xy}$ with $\dim_K C_{xy}\leq 2$ for all $x,y\in Q_0$.
Nevertheless, this fact does not imply that the quiver is intervally
finite (the number of paths between two vertices is finite). It is
enough to consider the quiver
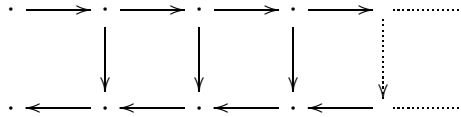
\[ \xymatrix { \cdot \ar[r] &\cdot\ar[r]\ar[d] & \cdot\ar[r]\ar[d]
&\cdot\ar[r]\ar[d] &\ar@{.>}[d] \ar@{.}[r] &  \\
\cdot &\cdot\ar[l] & \cdot \ar[l]&\cdot \ar[l]& \ar[l]\ar@{.}[r] &
}\] and the admissible coalgebra $C=C(Q,\Omega)$, where $\Omega$ is
the ideal
$$\Omega=KQ_2\oplus KQ_3\oplus \cdots \oplus KQ_n\oplus\cdots.$$ $C$
is a string coalgebra and then it is tame (see \cite[Section
6]{simson2}). \end{rem}

\end{document}